\def\yes{\if00}

\def\IfPic{\yes}

\documentstyle[leqno,12pt,amscd]{amsart}
\let\cal\mathcal
\theoremstyle{plain}
\newtheorem{Theorem}{Theorem}[section]
\newtheorem{Proposition}[Theorem]{Proposition}
\newtheorem{Lemma}[Theorem]{Lemma}
\newtheorem{Corollary}[Theorem]{Corollary}

\theoremstyle{definition}

\newtheorem{Remark}[Theorem]{Remark}

\newtheorem{Conjecture}[Theorem]{Conjecture}

\renewcommand{\theTheorem}{\arabic{section}.\arabic{Theorem}}
\renewcommand{\theClaim}{\arabic{section}.\arabic{Theorem}.\arabic{Claim}}
\renewcommand{\theequation}{\arabic{section}.\arabic{Theorem}.\arabic{Claim}}

\newcommand{\ZZ}{{\Bbb{Z}}}
\newcommand{\QQ}{{\Bbb{Q}}}
\newcommand{\RR}{{\Bbb{R}}}
\newcommand{\CC}{{\Bbb{C}}}
\newcommand{\PP}{{\Bbb{P}}}
\newcommand{\OO}{{\cal{O}}}
\newcommand{\Div}{\operatorname{Div}}
\newcommand{\CV}{\operatorname{CV}}
\newcommand{\Pic}{\operatorname{Pic}}
\newcommand{\Spec}{\operatorname{Spec}}
\newcommand{\rank}{\operatorname{rk}}

\newcommand{\Proof}{{\sl Proof.}\quad}
\newcommand{\QED}{{\unskip\nobreak\hfil\penalty50\quad\null\nobreak\hfil
{$\Box$}\parfillskip0pt\finalhyphendemerits0\par\medskip}}
\newcommand{\rest}[2]{\left.{#1}\right\vert_{{#2}}}

\IfPic
\input DraTex.sty
\input AlDraTex.sty
\else
\fi

\begin{document}

\title[Bogomolov conjecture over function fields]
{Bogomolov conjecture over function fields \\
for stable curves with only irreducible fibers}
\author{Atsushi Moriwaki}
\date{June, 1995 (Version 2.0)}
\address{Department of Mathematics, Faculty of Science,
Kyoto University, Kyoto, 606-01, Japan}
\curraddr{Department of Mathematics, University of California,
Los Angeles, 405 Hilgard Avenue, Los Angeles, CA 90024, USA}
\email{moriwaki@@math.ucla.edu}
\thanks{We would like to thank Prof. S. Zhang for his helpful comments.}
\begin{abstract}
Let $K$ be a function field and $C$ a non-isotrivial curve of genus $g \geq 2$
over
$K$. In this paper, we will show that if $C$ has a global stable model
with only geometrically irreducible fibers,
then Bogomolov conjecture over function fields holds.
\end{abstract}

\maketitle

\tableofcontents

\section{Introduction}

Let $k$ be a field,
$X$ a smooth projective surface over $k$, $Y$ a smooth projective curve over
$k$,
and $f : X \to Y$ a generically smooth semistable curve of genus $g \geq 2$
over $Y$.
Let $K$ be the function field of $Y$, $\overline{K}$
the algebraic closure of $K$, and $C$ the generic fiber of $f$.
For $D \in \Pic^1(C)(\overline{K})$, let
$j_D : C_{\overline{K}} \to \Pic^0(C)_{\overline{K}}$
be an embedding defined by $j_D(x) = x - D$.
Then, we have the following conjecture due to Bogomolov.

\begin{Conjecture}[Bogomolov conjecture over function fields]
\label{conj:Geometric:Bogomolov:Conjecture}
If $f$ is non-isotrivial, then, for any
embedding $j_D$, the image $j_D(C(\overline{K}))$ is discrete in terms
of the semi-norm $\Vert \ \Vert_{NT}$ given by
the Neron-Tate height pairing on $\Pic^0(C)(\overline{K})$, i.e.,
for any point $P \in \Pic^0(C)(\overline{K})$,
there is a positive number $\epsilon$ such that the set
\[
 \left\{ x \in C(\overline{K}) \mid
         \Vert j_D(x) - P \Vert_{NT} \leq \epsilon \right\}
\]
is finite.
\end{Conjecture}

In this paper, we will prove the above conjecture
under the assumption that
the stable model of
$f : X \to Y$ has only geometrically irreducible fibers.
\begin{Theorem}
\label{thm:conj:bogomolov}
If the stable model of $f : X \to Y$ has only geometrically irreducible fibers,
then Conjecture~\ref{conj:Geometric:Bogomolov:Conjecture} holds.
More strongly, there is  a positive number $A$ with the following properties.
\begin{enumerate}
\renewcommand{\labelenumi}{(\arabic{enumi})}
\item
${\displaystyle
A \geq \sqrt{\frac{g-1}{12g(2g+1)}\delta}}$,
where $\delta$ is the number of singularities in singular fibers
of $f_{\bar{k}} : X_{\bar{k}} \to Y_{\bar{k}}$.

\item
For any small positive number $\epsilon$, the set
\[
 \left\{ x \in C(\overline{K}) \mid \Vert j_D(x) - P \Vert_{NT} \leq
 (1 - \epsilon)A \right\}
\]
is finite for any embedding $j_D$ and
any point $P \in \Pic^0(C)(\overline{K})$.
\end{enumerate}
\end{Theorem}

Our proof of Theorem~\ref{thm:conj:bogomolov}
is based on the admissible pairing on semistable curves due to S. Zhang (cf.
\S\ref{sec:metrized:graph:green:function:admissible:pairing}),
Cornalba-Harris-Xiao's inequality over an arbitrary field
(cf. Theorem~\ref{thm:Cornalba-Harris-Xiao:inequality}) and
an exact calculation of a Green function on
a certain metrized graph (cf. Lemma~\ref{lem:green:union:circle}).
The estimation of a Green function also gives the following
result, which strengthen S. Zhang's theorem \cite{Zh}.

\begin{Theorem}[cf. Corollary~\ref{cor:lower:bound:w:w:not:smooth}]
Let $K$ be a number field, $O_K$ the ring of integers,
$f : X \to \Spec(O_K)$ a regular semistable arithmetic surface of genus $g \geq
2$
over $O_K$.
If $f$ is not smooth, then
\[
(\omega_{X/O_K}^{Ar} \cdot \omega_{X/O_K}^{Ar}) \geq \frac{\log 2}{6(g-1)}.
\]
\end{Theorem}

\section{Metrized graph, Green function and admissible pairing}
\label{sec:metrized:graph:green:function:admissible:pairing}

In this section, we recall several facts of metrized graphs,
Green functions and the admissible pairing
on semistable curves. Details can be found in Zhang's paper \cite{Zh}.

Let $G$ be a locally metrized and compact topological space.
We say $G$ is a metrized graph if, for any $x \in G$,
there is a positive number $\epsilon$, a positive integer $d = v(x)$
(which is called the valence at $x$), and
an open neighborhood $U$ of $x$ such that $U$ is isometric to
\[
  \left\{ t e^{\frac{2 \pi \sqrt{-1} k}{d}} \in \CC \mid 0 \leq t < \epsilon,
  k \in \ZZ \right\}.
\]
Let $\Div(G)$ be a free abelian group generated by points of $G$.
An element of $\Div(G)$ is called a {\em divisor} on $G$.
Let $F(G)$ be the set of all piecewisely smooth real valued
functions on $G$.
For $f \in F(G)$, we can define the Dirac function $\delta(f)$ associated with
$f$
as follows. If $x \in G$ and $v(x) = n$, then $\delta(f)(x)$ is given by
\[
   (\delta(f)(x), g) = g(x) \sum_{i=1}^n \lim_{x_i \to 0} f'(x_i),
\]
where $g \in F(G)$ and $x_i$ is the arc-length parameter of one branch
starting from $x$. The Laplacian $\Delta$ for $f \in F(G)$ is defined by
\[
\Delta(f) = -f'' - \delta(f),
\]
where $f''$ is the second derivative of $f$ in the sense of distribution.
Let $Q(G)$ be a subset of $F(G)$ consisting of piecewisely
quadric polynomial functions.
Let $V$ be a set of vertices of $G$ such that
$G \setminus V$ is a disjoint union of open segments.
Let $E$ be the collection of segments in $G \setminus V$.
We denote by $Q(G, V)$ a subspace of $Q(G)$ consisting of
functions whose restriction to each edge in $E$ are
quadric polynomial functions, and by $M(G, V)$
a vector space of measures on $G$ generated by Dirac functions $\delta_v$ at
$v \in V$ and by Lebesgue measures on edges $e \in E$ arising from
the arc-length parameter.
The fundamental theorem is the following existence of
the admissible metric and the Green function.

\begin{Theorem}[{\cite[Theorem 3.2]{Zh}}]
\label{thm:existence:metric:green}
Let $D = \sum_{x \in G} d_x x$ be a divisor on $G$ such that
the support of $D$ is in $V$.
If $G$ is connected and $\deg(D) \not= -2$, then
there are a unique measure $\mu \in M(G, V)$ and
a unique function $g_{\mu}$ on $G \times G$ with the following properties.
\begin{enumerate}
\renewcommand{\labelenumi}{(\arabic{enumi})}
\item
${\displaystyle \int_{G} \mu = 1}$.

\item
$g_{\mu}(x, y)$ is symmetric and continuous on $G \times G$.

\item
For a fixed $x \in G$, $g_{\mu}(x, y) \in Q(G)$.
Moreover, if $x \in V$, then $g_{\mu}(x, y) \in Q(G, V)$.

\item
For a fixed $x \in G$, $\Delta_y(g_{\mu}(x, y)) = \delta_x - \mu$.

\item
For a fixed $x \in G$, ${\displaystyle \int_G g_{\mu} (x, y) \mu(y) = 0}$.

\item
$g_{\mu}(D, y) + g_{\mu}(y, y)$ is a constant for all $y \in G$,
where $g_{\mu}(D, y) = \sum_{x \in G} d_x g_{\mu}(x, y)$.
\end{enumerate}
Further, if $d_x \geq v(x) - 2$ for all $x \in G$,
then $\mu$ is positive.
\end{Theorem}

The measure $\mu$ in Theorem~\ref{thm:existence:metric:green} is called
the {\em admissible metric} with respect to $D$ and $g_{\mu}$ is called
the {\em Green function} with respect to $\mu$.
The constant $g_{\mu}(D, y) + g_{\mu}(y, y)$ is denoted by $c(G, D)$.

\bigskip
Let $k$ be an algebraically closed field,
$X$ a smooth projective surface over $k$, $Y$ a smooth projective curve over
$k$,
and $f : X \to Y$ a generically smooth semistable curve of genus $g \geq 1$
over $Y$.
Let $\CV(f)$ be the set of all critical values of $f$, i.e.,
$y \in \CV(f)$ if and only if $f^{-1}(y)$ is singular.
For $y \in \CV(f)$, let $G_y$ be the metrized graph of $f^{-1}(y)$
defined as follows.
The set of vertices $V_y$ of $G_y$ is indexed by
irreducible components of the fiber $f^{-1}(y)$ and
singularities of $f^{-1}(y)$ correspond to edges of length $1$.
We denote by $C_v$ the corresponding irreducible curve for a vertex $v$
in $V_y$. Let $K_y$ be a divisor on $G_y$ given by
\[
  K_y = \sum_{v \in V_y} (\omega_{X/Y} \cdot C_v) v.
\]
Let $\mu_y$ be the admissible metric with respect to $K_y$ and
$g_{\mu_y}$ the Green function of $\mu_y$.
The admissible dualizing sheaf $\omega_{X/Y}^a$ is defined by
\[
\omega_{X/Y}^a = \omega_{X/Y} - \sum_{y \in \CV(f)} c(G_y, K_y) f^{-1}(y).
\]
Here we define a new pairing $(D \cdot E)_a$ for
$D, E \in \Div(X) \otimes \RR$ by
\[
(D \cdot E)_a = (D \cdot E) + \sum_{y \in \CV(f)} \left\{ \sum_{v, v' \in V_y}
(D \cdot C_v) g_{\mu_y}(v, v') (E \cdot C_{v'}) \right\}.
\]
This pairing is called the {\em admissible pairing}.
It has lots of properties. For our purpose, the following are important.

\medskip
\begin{enumerate}
\renewcommand{\labelenumi}{(\arabic{enumi})}
\item (Adjunction formula)
If $B$ is a section of $f$, then
$(\omega_{X/Y}^a + B \cdot B)_a = 0$.

\item (Intersection with a fiber)
If $D$ is an $\RR$-divisor with degree $0$ along general fibers,
then $(D \cdot Z)_a = 0$ for all vertical curves $Z$.
(cf. Proposition~\ref{prop:admissible:with:fiber})

\item (Compatibility with base changes)
The admissible pairing is compatible with base changes.
Namely, let $\pi : Y' \to Y$ be a finite morphism of smooth projective curves,
and
$X'$ the minimal resolution of the fiber product of $X \times_Y Y'$.
We set the induced morphisms as follows.
\[
\begin{CD}
X @<{\pi'}<< X' \\
@V{f}VV  @VV{f'}V \\
Y @<<{\pi}<  Y'
\end{CD}
\]
Then, for $D, E \in \Div(X) \otimes \RR$,
$({\pi'}^*(D) \cdot {\pi'}^*(E))_a = (\deg \pi)(D \cdot E)_a$.
Moreover, we have ${\pi'}^*(\omega_{X/Y}^a) = \omega_{X'/Y'}^a$.
Thus, $(\omega_{X'/Y'}^a \cdot \omega_{X'/Y'}^a)_a
= (\deg \pi) (\omega_{X/Y}^a \cdot \omega_{X/Y}^a)_a$.
\end{enumerate}

\medskip\noindent
Using the above properties, we can give the Neron-Tate height paring
in terms of the admissible pairing.
Let $C$ be the generic fiber of $f$, $K$ the function field of $Y$, and
$L, M \in \Pic^0(C)(\overline{K})$. Then, there are
a base change $Y' \to Y$, a semistable model $X'$ of $C$ over $Y'$, and
line bundles ${\cal L}$ and ${\cal M}$ on $X'$ such that
${\cal L}_{\overline{K}} = L$ and ${\cal M}_{\overline{K}} = M$.
Moreover, we can find vertical $\QQ$-divisors $V$ and $V'$ on $X'$
such that $({\cal L} + V \cdot Z) = ({\cal M} + V' \cdot Z) = 0$
for all vertical curves $Z$ on $X'$. Then, it is easy to see that
\[
  \frac{-1}{[k(Y') : k(Y)]}({\cal L} + V \cdot {\cal M} + V')
\]
is well-defined. It is denoted by $(L \cdot M)_{NT}$ and is called
the {\em Neron Tate height pairing}. Moreover, it is easy to see
$(L \cdot L)_{NT} \geq 0$. So $\sqrt{(L \cdot L)_{NT}}$ is denoted
by $\Vert L \Vert_{NT}$. On the other hand, by the definition of the
admissible pairing, we have
\[
({\cal L} + V \cdot {\cal M} + V') = ({\cal L} + V \cdot {\cal M} + V')_a.
\]
Thus, using the second property of the above, we can see that
\[
-[k(Y') : k(Y)] (L \cdot M)_{NT} = ({\cal L}\cdot {\cal M})_a,
\]
which means that the admissible pairing does not depend
on the choice of the compactification of $L$ and $M$, and that
of course
\[
(L \cdot M)_{NT} = \frac{-({\cal L}\cdot {\cal M})_a}{[k(Y') : k(Y)]}.
\]

\medskip
Next, let us consider a height function in terms of the
admissible pairing. Let ${\cal L}$ be an $\RR$-divisor on $X$ and
$x \in C(\overline{K})$. Then, taking a suitable base change $\pi : Y' \to Y$,
there is a semistable model
$f' : X' \to Y'$ of $C$ such that $x$ is realized as a section $B_x$
of $f'$. We set
\[
   h^a_{{\cal L}}(x) = \frac{({\pi'}^*({\cal L}) \cdot B_x)_a}{\deg \pi},
\]
where $\pi' : X' \to X$ is the induced morphism.
One can easily see that $h^a_{{\cal L}}(x)$ is well-defined by the third
property
of the above. The following generic lower estimate of the height function
is important for our purpose.

\begin{Theorem}[{\cite[Theorem 5.3]{Zh}}]
\label{thm:lower:estimate:height}
If $\deg({\cal L}_K) > 0$ and ${\cal L}$ is $f$-nef, then,
for any $\epsilon > 0$, there is a finite subset $S$ of $C(\overline{K})$
such that
\[
h^a_{{\cal L}}(x) \geq \frac{({\cal L} \cdot {\cal L})_a}{2 \deg({\cal L}_K)} -
\epsilon
\]
for all $x \in C(\overline{K}) \setminus S$.
\end{Theorem}

As corollary, we have the following.

\begin{Corollary}[{\cite[Theorem 5.6]{Zh}}]
\label{cor:lower:estimate:NT:metric}
Let $D \in \Pic^1(C)(\overline{K})$. Then, for any $\epsilon > 0$,
there is a finite subset $S$ of $C(\overline{K})$
such that
\[
\Vert D - x \Vert_{NT}^2 \geq
\frac{(\omega_{X/Y}^a \cdot \omega_{X/Y}^a)_a}{4(g-1)} +
\frac{\Vert \omega_C - (2g-2) D \Vert_{NT}^2}{4g(g-1)}
- \epsilon
\]
for all $x \in C(\overline{K}) \setminus S$.
\end{Corollary}

\Proof
Let $\pi_1 : Y_1 \to Y$ be a base change of $f : X \to Y$ such that
$D$ is defined over the function field $k(Y_1)$ of $Y_1$.
Let $f_1 : X_1 \to Y_1$ be the semistable model of $C$ over $Y_1$,
$F$ a general fiber of $f_1$, and
${\cal D}$ a compactification of $D$ such that
${\cal D}$ is a $\QQ$-divisor on $X_1$ and ${\cal D}$ is $f_1$-nef.
Using adjuction formula and applying Theorem~\ref{thm:lower:estimate:height} to
\[
  {\cal L} = \omega_{X_1/Y_1}^a + 2 {\cal D} - ({\cal D} \cdot {\cal D})_a F,
\]
we have our corollary.
\QED

\section{Green function of a certain metrized graph}
\label{sec:certain:graph}

In this section, we will construct a Green function of
a certain metrized graph. Let us begin with the following lemma.

\begin{Lemma}
\label{lem:laplacian:on:circle}
Let $C$ be a circle with arc-length $l$.
Fixing a point $O$ on $C$, let $t : C \to [0, l)$ be a coordinate of $C$
with $t(O) = 0$ coming from an arc-length parameterization of $C$.
We set
\[
\phi(t) = \frac{1}{2l} t^2 - \frac{1}{2} |t|
\quad\text{and}\quad
f(x, y) = \phi(t(x) - t(y)).
\]
Then, we have the following.
\begin{enumerate}
\renewcommand{\labelenumi}{(\arabic{enumi})}
\item
$f(x, y)$ is symmetric and continuous on $C \times C$.

\item
$f(x, y)$ is smooth on the outside of the diagonal.

\item
For a fixed $x \in C$,
${\displaystyle \Delta_y(f(x, y)) = \delta_x - \frac{dt}{l}}$.
\end{enumerate}
\end{Lemma}

\Proof
We can check them by a straightforward calculation.
\QED

Let $C_1, \ldots, C_n$ be circles and $G$ a metrized graph constructed by
joining $C_i$'s at a point $O$. Let $l_i$ be the arc-length of $C_i$ and
$t_i : C_i \to [0, l_i)$ a coordinate of $C_i$ with $t_i(O) = 0$.
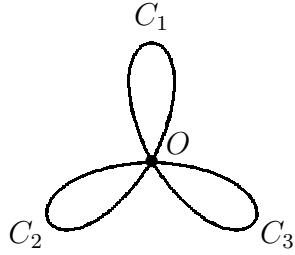
\IfPic
\par\bigskip
\Draw
\MoveTo(0,0) \MarkLoc(O) \Node(Q)(--$\bullet$--)
\MoveTo(30,60) \MarkLoc(A)
\MoveTo(-30,60) \MarkLoc(B)
\MoveTo(-66.96,-4.02) \MarkLoc(C)
\MoveTo(-36.96,-55.98) \MarkLoc(D)
\MoveTo(36.96,-55.98) \MarkLoc(E)
\MoveTo(66.96,-4.02) \MarkLoc(F)
\Curve(O,A,B,O) \Curve(O,C,D,O) \Curve(O,E,F,O)
\MoveTo(10,7) \Node(G)(--$O$--)
\MoveTo(0,55) \Node(H)(--$C_1$--)
\MoveTo(-47.63,-27.5) \Node(H)(--$C_2$--)
\MoveTo(47.63,-27.5) \Node(H)(--$C_3$--)
\EndDraw
\bigskip
\par\noindent
\else
\fi
{}From now on, we will identify a point on $C_i$ with its coordinate.
As in Lemma~\ref{lem:laplacian:on:circle}, for each $i$, we set
\[
  \phi_i(t) = \frac{1}{2l_i} t^2 - \frac{1}{2} |t|.
\]
We fix a positive integer $g$.
Here we consider a measure $\mu$ and a divisor $K$ on $G$ defined by
\[
    \mu = \frac{g-n}{g} \delta_O + \sum_{i=1}^n \frac{d t_i}{gl_i}
    \quad\text{and}\quad
    K = (2g-2)O.
\]
Moreover, let us consider
the following function $g_{\mu}$ on $G \times G$.
\[
g_{\mu}(x, y) =
\begin{cases}
{\displaystyle \phi_i(x - y) - \frac{g-1}{g} \left(\phi_i(x) + \phi_i(y)\right)
+
\frac{L}{12g^2}} &
\text{if $x, y \in C_i$} \\
{\displaystyle \frac{1}{g}\left( \phi_i(x) + \phi_j(y) \right) +
\frac{L}{12g^2}} &
\text{if $x \in C_i$, $y \in C_j$ and $i \not= j$}
\end{cases}
\]
where $L = l_1 + \cdots + l_n$.
Then, we can see the following.

\begin{Lemma}
\label{lem:green:union:circle}
\begin{enumerate}
\renewcommand{\labelenumi}{(\arabic{enumi})}
\item ${\displaystyle \int_G \mu = 1}$.

\item
$g_{\mu}(x, y)$ is symmetric and continuous on $G \times G$.

\item
For a fixed $x \in G$, $\Delta_y(g_{\mu}(x, y)) = \delta_x - \mu$.

\item
For a fixed $x \in G$, ${\displaystyle \int_G g_{\mu}(x, y) \mu(y) = 0}$.

\item
${\displaystyle g_{\mu}(K, y) + g_{\mu}(y, y) = \frac{L(2g-1)}{12g^2}}$
for all $y \in G$.
\end{enumerate}
\end{Lemma}

\Proof
(1), (2) These are obvious.

(3) We assume $x \in C_i$. By \cite[Lemma a.4, (a)]{Zh},
\[
\Delta_y(g_{\mu}(x, y)) = \sum_{j=1}^n \Delta_y(\rest{g_{\mu}(x, y)}{C_j}).
\]
Therefore, using Lemma~\ref{lem:laplacian:on:circle},
we get
\begin{align*}
\Delta_y(g_{\mu}(x, y)) & =
\Delta_y(\rest{g_{\mu}(x, y)}{C_i}) + \sum_{j \not= i}^n
\Delta_y(\rest{g_{\mu}(x, y)}{C_j}) \\
& = \left( \delta_x - \frac{dt_i}{l_i} -
\frac{g-1}{g}\left(\delta_O - \frac{d t_i}{l_i}\right)\right) +
\sum_{j \not= i}^n \frac{1}{g}\left( \delta_O - \frac{dt_j}{l_j} \right) \\
& = \delta_x - \mu.
\end{align*}

(4) We assume $x \in C_i$.
Then, by a direct calculation, we can see
\[
\int_{C_j} g_{\mu}(x, t_j) \frac{d t_j}{gl_j} =
\begin{cases}
{\displaystyle -\frac{g-1}{g^2}\phi_i(x) - \frac{l_i}{12g^2} + \frac{L}{12g^3}}
&
\text{if $j = i$} \\
{\displaystyle \frac{1}{g^2} \phi_i(x) - \frac{l_j}{12g^2} + \frac{L}{12g^3}} &
\text{if $j \not= i$}
\end{cases}
\]
Therefore,
\[
\sum_{j=1}^n \int_{C_j} g_{\mu}(x, t_j) \frac{d t_j}{gl_j} =
\frac{n-g}{g}\left( \frac{1}{g} \phi_i(x) + \frac{L}{12g^2} \right).
\]
Hence,
\begin{align*}
\int_G g_{\mu}(x, y) \mu(y) & =
\frac{g-n}{g}g_{\mu}(x, 0) + \sum_{j=1}^n \int_{C_j} g_{\mu}(x, t_j) \frac{d
t_j}{gl_j} \\
& = \frac{g-n}{g}\left( \frac{1}{g} \phi_i(x) + \frac{L}{12g^2} \right) +
\frac{n-g}{g}\left( \frac{1}{g} \phi_i(x) + \frac{L}{12g^2} \right) \\
& = 0.
\end{align*}

(5) Since
\[
 g_{\mu}(O, x) = \frac{1}{g} \phi_i(x) + \frac{L}{12g^2} \quad\text{and}\quad
 g_{\mu}(x, x) = \frac{-2(g-1)}{g} \phi_i(x) + \frac{L}{12g^2},
\]
(5) follows.
\QED

This lemma says us that
$\mu$ is the admissible metric with respect to $K$,
$g_{\mu}$ is the Green function of $\mu$, and
${\displaystyle c(G, K) = \frac{L(2g-1)}{12g^2}}$.

\section{Cornalba-Harris-Xiao's inequality over an arbitrary field}

In this section, we would like to generalize Cornalba-Harris-Xiao's inequality
to fibered algebraic surfaces over an arbitrary field, namely,

\begin{Theorem}
\label{thm:Cornalba-Harris-Xiao:inequality}
Let $k$ be a field,
$X$ a smooth projective surface over $k$, $Y$ a smooth projective curve over
$k$, and
$f : X \to Y$ a generically smooth morphism with $f_*\OO_X = \OO_C$.
If the genus $g$ of the generic fiber of $f$ is greater than or equal to $2$
and
$\omega_{X/Y}$ is $f$-nef, then
\[
(\omega_{X/Y} \cdot \omega_{X/Y}) \geq \frac{4(g-1)}{g}
\deg(f_*(\omega_{X/Y})).
\]
\end{Theorem}

The above was proved in \cite{CH} and \cite{Xi} under the assumption
$\operatorname{char}(k) = 0$. Here we prove it using the following result
of Bost.

\begin{Theorem}[{\cite[Theorem III]{Bo}}]
\label{thm:bost:ineq}
Let $k$ be a field, $Y$ a smooth projective curve over $k$, and
$E$ a vector bundle on $Y$. Let
\[
\pi : P =
\operatorname{Proj}\left(
\bigoplus_{n=0}^{\infty} \operatorname{Sym}^n(E)
\right)
\longrightarrow Y
\]
be the projective bundle of $E$ and $\OO_P(1)$ the tautological line bundle
on $P$.
If an effective cycle $Z$ of dimension $d \geq 1$ on $P$ is
Chow semistable on the generic fiber of $\pi$, then
\[
\frac{\left(\OO_P(1)^d \cdot Z \right)}
{d \cdot \left( \OO_P(1)^{d-1} \cdot Z \cdot F \right)}
\geq \frac{\deg E}{\rank E},
\]
where $F$ is a general fiber of $\pi$.
\end{Theorem}

\medskip
First of all, let us begin with the following lemmas.

\begin{Lemma}
\label{lem:Chow:semistable:canonical}
Let $K$ be a field, $C$ a smooth projective curve over $K$ of genus $g \geq 2$,
and
$\phi : C \to \PP^{g -1}$ a morphism given by
the complete linear system $|\omega_C|$. Then $\phi_*(C)$
is a Chow semistable cycle on $\PP^{g-1}$.
\end{Lemma}

\proof
Let $R$ be the image of $C$ by $\phi$ and $n$ an integer given by
\[
  n =
  \begin{cases}
    1, & \text{if $C$ is non-hyperelliptic}, \\
    2, & \text{if $C$ is hyperelliptic}.
  \end{cases}
\]
Then, $\phi_*(C) = nR$.
Thus a Chow form of $\phi_*(C)$ is the $n$-th power of
a Chow form of $R$. Therefore, $\phi_*(C)$ is Chow semistable if and only if
$R$ is Chow semistable.
Moreover, Theorem~4.12 in \cite{Mu} says that
Chow semistability of $R$ is derived from linear semistability of $R$ .

Let $V$ be a subspace of $H^0(C, \omega_C)$,
$p : \PP^{g-1} \dashrightarrow \PP^{\dim V - 1}$
the projection defined by the inclusion $V \hookrightarrow H^0(C, \omega_C)$,
and $\phi' : C \to \PP^{\dim V - 1}$ a morphism given by $V$.
Then, $p \cdot \phi = \phi'$.
We need to show that
\addtocounter{Claim}{1}
\begin{equation}
\label{eqn:lem:stable:kernel:Chow:stable}
\frac{2}{n} = \frac{\deg(R)}{g - 1} \leq \frac{\deg(p_*(R))}{\dim V - 1}
\end{equation}
to see linear semistability of $R$.
Since $\deg({\phi'}^*(\OO(1))) = n \deg(p_*(R))$,
(\ref{eqn:lem:stable:kernel:Chow:stable})
is equivalent to say
\[
2 \leq
\frac{\deg({\phi'}^*(\OO(1)))}{\dim V - 1}.
\]
On the other hand,
if we denote by
$\omega^V_C$ the image of $V \otimes \OO_C \to \omega_C$,
then, by Clifford's lemma, we have
\[
\dim V - 1 \leq \dim |\omega^V_C| \leq \frac{\deg(\omega^V_C)}{2}.
\]
Thus, we get (\ref{eqn:lem:stable:kernel:Chow:stable})
because
$\deg({\phi'}^*(\OO(1))) = \deg(\omega^V_C)$.
\QED

\begin{Remark}
By \cite[Proposition~4.2]{Bo}, $\phi_*(C)$ is
actually Chow stable when $\operatorname{char}(K) = 0$.
We don't know whether $\phi_*(C)$ is Chow stable if $\operatorname{char}(K) >
0$.
Anyway, semistability is enough for our purpose.
\end{Remark}

\medskip
Let us start the proof of Theorem~\ref{thm:Cornalba-Harris-Xiao:inequality}.
Let
\[
\phi : X \dashrightarrow P = \operatorname{Proj}\left(
\bigoplus_{n=0}^{\infty} \operatorname{Sym}^n(f_*(\omega_{X/Y}))
\right)
\]
be a rational map over $Y$ induced by $f^*f_*(\omega_{X/Y}) \to \omega_{X/Y}$.
Here we take a birational morphism $\mu : X' \to X$ of smooth projective
varieties such that $\phi' = \phi \cdot \mu : X' \to P$
is a morphism.
\[
\begin{CD}
X' @>{\mu}>> X \\
@V{\phi'}VV @VV{\phi}V \\
P @= P
\end{CD}
\]
Then, there is an effective vertical divisor $D$ on $X'$ such that
$\mu^*(\omega_{X/Y}) = {\phi'}^*(\OO_P(1)) +  D$.
Let $Z = {\phi'}_*(X')$. Then, by Lemma~\ref{lem:Chow:semistable:canonical},
$Z$ give a Chow semistable cycle on the generic fiber.
Thus, by Theorem~\ref{thm:bost:ineq}, we have
\[
\frac{\left({\phi'}^*(\OO_P(1)) \cdot {\phi'}^*(\OO_P(1))\right)}{4(g-1)}
\geq \frac{\deg(f_*(\omega_{X/Y}))}{g}.
\]
On the other hand, since $\omega_{X/Y}$ is $f$-nef and $(D \cdot D) \leq 0$,
\begin{align*}
\left({\phi'}^*(\OO_P(1)) \cdot {\phi'}^*(\OO_P(1))\right) & =
\left( \mu^*(\omega_{X/Y}) - D \cdot \mu^*(\omega_{X/Y}) - D \right) \\
& = \left( \omega_{X/Y} \cdot \omega_{X/Y} \right)
-2 \left(\mu^*(\omega_{X/Y}) \cdot D \right) + (D \cdot D) \\
& \leq \left( \omega_{X/Y} \cdot \omega_{X/Y} \right).
\end{align*}
Therefore, we have our desired inequality.
\QED

\begin{Remark}
\label{rem:semistable:kernel:another:proof}
If $\operatorname{char}(k) = 0$,
we can give another proof of Theorem~\ref{thm:Cornalba-Harris-Xiao:inequality}
according to \cite{Mo2}. A rough idea is the following.
Since the kernel $K$ of $f^*f_*(\omega_{X/Y}) \to \omega_{X/Y}$
is semistable on the generic fiber of $f$ by virtue of \cite{PR}, we can apply
Bogomolov-Gieseker's inequality to $K$, which implies
Cornalba-Harris-Xiao's inequality by easy calculations.
\end{Remark}

\section{Proof of Theorem~\ref{thm:conj:bogomolov}}

In this section, we would like to give the proof of
Theorem~\ref{thm:conj:bogomolov}.
First of all, let us fix notations.
Let $k$ be a field,
$X$ a smooth projective surface over $k$, $Y$ a smooth projective curve over
$k$,
and $f : X \to Y$ a generically smooth semistable curve of genus $g \geq 2$
over $Y$.
Let $K$ be the function field of $Y$,
$\overline{K}$ the algebraic closure of $K$, and $C$ the generic fiber of $f$.
We assume that $f$ is non-isotrivial and
the stable model of $f : X \to Y$ has only geometrically irreducible fibers.
Clearly, for the proof of Theorem~\ref{thm:conj:bogomolov},
we may assume that $k$ is algebraically closed.
Then, we have the following lower estimate of
$(\omega_{X/Y}^a \cdot \omega_{X/Y}^a)_a$.

\begin{Theorem}
\label{thm:lower:bound:admissible:intersection}
Under the above assumptions, $(\omega_{X/Y}^a \cdot \omega_{X/Y}^a)_a$
is positive. Moreover,
\[
 (\omega_{X/Y}^a \cdot \omega_{X/Y}^a)_a \geq \frac{(g-1)^2}{3g(2g+1)} \delta,
\]
where $\delta$ is the number of singularities in singular fibers of $f$.
\end{Theorem}

\Proof
Let $\CV(f)$ be the set of all critical values of $f$.
For $y \in \CV(f)$, the number of singularities of $f^{-1}(y)$
is denoted by $\delta_y$.
Let $G_y$ be the metrized graph of $f^{-1}(y)$ as in
\S\ref{sec:metrized:graph:green:function:admissible:pairing}.
Then, the total arc-length of $G_y$ is $\delta_y$.
Let $K_y$ be the divisor on $G_y$ coming from $\omega_{X/Y}$
as in \S\ref{sec:metrized:graph:green:function:admissible:pairing},
$\mu_y$ the admissible metric of $K_y$, and $g_{\mu_y}$ the Green
function of $\mu_y$. By the definition of $\omega_{X/Y}^a$
(see \S\ref{sec:metrized:graph:green:function:admissible:pairing}), we have
\[
(\omega_{X/Y}^a \cdot \omega_{X/Y}^a)_a =
(\omega_{X/Y} \cdot \omega_{X/Y}) +
\sum_{ y \in \CV(f)} \left\{ g_{\mu_y}(K_y, K_y) - 2(2g-2) c(G_y, K_y)
\right\}.
\]
On the other hand, $G_y$ is isometric to the graph treated in
\S\ref{sec:certain:graph}. Thus, by Lemma~\ref{lem:green:union:circle},
\begin{align*}
g_{\mu_y}(K_y, K_y) - 2(2g-2) c(G_y, K_y) & =
(2g-2)^2 \frac{\delta_y}{12g^2} - 2(2g-2) \frac{(2g-1) \delta_y}{12g^2} \\
& = -\frac{g-1}{3g} \delta_y.
\end{align*}
Thus
\[
(\omega_{X/Y}^a \cdot \omega_{X/Y}^a)_a =
(\omega_{X/Y} \cdot \omega_{X/Y}) - \frac{g-1}{3g} \delta.
\]
By virtue of Theorem~\ref{thm:Cornalba-Harris-Xiao:inequality} and
Noether formula
\[
\deg(f_*(\omega_{X/Y})) =
\frac{(\omega_{X/Y} \cdot \omega_{X/Y}) + \delta}{12},
\]
we have
\[
  ( \omega_{X/Y} \cdot \omega_{X/Y} ) \geq \frac{g-1}{2g+1} \delta.
\]
Therefore, we get
\[
 (\omega_{X/Y}^a \cdot \omega_{X/Y}^a)_a \geq \frac{(g-1)^2}{3g(2g+1)} \delta.
\]
In particular, $(\omega_{X/Y}^a \cdot \omega_{X/Y}^a)_a > 0$
if $f$ is not smooth. Further, if $f$ is smooth, then
\[
(\omega_{X/Y}^a \cdot \omega_{X/Y}^a)_a = (\omega_{X/Y} \cdot \omega_{X/Y}) > 0
\]
because $f$ is non-isotrivial.
\QED

\bigskip
Let us start the proof of Theorem~\ref{thm:conj:bogomolov}.
We set
\[
    A = \sqrt{\frac{(\omega_{X/Y}^a \cdot \omega_{X/Y}^a)_a}{4(g-1)}}.
\]
Then, by Theorem~\ref{thm:lower:bound:admissible:intersection},
$A$ is positive and
\[
A \geq \sqrt{\frac{g-1}{12g(2g+1)} \delta}.
\]
By virtue of Corollary~\ref{cor:lower:estimate:NT:metric},
for any $D \in \Pic^1(C)(\overline{K})$ and
any $P \in \Pic^0(C)(\overline{K})$, there is a finite subset $S$ of
$C(\overline{K})$ such that
\[
\Vert x - D - P \Vert_{NT} > (1-\epsilon)A
\]
for all $x \in C(\overline{K}) \setminus S$.
Therefore, we have
\[
 \left\{ x \in C(\overline{K}) \mid \Vert j_D(x) - P \Vert_{NT} \leq
 (1 - \epsilon)A \right\} \subset S.
\]
Thus, we get the second property of $A$.
\QED

\section{Effective lower bound of $(\omega \cdot \omega)$
for arithmetic surfaces}

Let $K$ be a number field, $O_K$ the ring of integers,
$f : X \to \Spec(O_K)$ a regular semistable arithmetic surface of genus $g \geq
2$
over $O_K$. In \cite{Mo3}, we proved the following.

\begin{Theorem}
If geometric fibers $X_{\overline{P}_1}, \ldots, X_{\overline{P}_n}$ of $X$
at $P_1, \ldots, P_n \in \Spec(O_K)$ are reducible, then
\[
\left(\omega_{X/O_K}^{Ar} \cdot \omega_{X/O_K}^{Ar} \right) \geq
\sum_{i=1}^n \frac{\log \#(O_K/P_i) }{6(g-1)}.
\]
\end{Theorem}

Using Lemma~\ref{lem:green:union:circle},
we have the following exact lower estimate for
stable curves with only irreducible fibers.

\begin{Theorem}
Assume that the stable model of $f : X \to \Spec(O_K)$ has only
geometric irreducible fibers.
If $\{ P_1, \ldots, P_n \}$ is the set of critical values of $f$, then
\[
(\omega_{X/O_K}^{Ar} \cdot \omega_{X/O_K}^{Ar})
\geq \sum_{i=1}^n \frac{g-1}{3g} \delta_i \log\#(O_K/P_i),
\]
where $\delta_i$ is the number of singularities of the geometric fiber
at $P_i$. Moreover, equality holds if and only if
there is a sequence of distinct points $x_1, x_2, \ldots$ of
$X(\overline{\QQ})$
such that
\[
\lim_{i \to \infty} \Vert (2g-2) x_i - \omega \Vert_{NT} = 0.
\]
\end{Theorem}

\Proof
By virtue of Lemma~\ref{lem:green:union:circle},
\[
(\omega_{X/O_K}^a \cdot \omega_{X/O_K}^a)_a =
(\omega_{X/O_K}^{Ar} \cdot \omega_{X/O_K}^{Ar}) -
\sum_{i=1}^n \frac{g-1}{3g} \delta_i \log \#(O_K/P_i).
\]
Therefore, our theorem follows from \cite[Corollary 5.7]{Zh}.
\QED

Combining the above two theorems, we have the following corollary,
which is a stronger version of S. Zhang's result \cite{Zh}.

\begin{Corollary}
\label{cor:lower:bound:w:w:not:smooth}
If $f : X \to \Spec(O_K)$ is not smooth, then
\[
(\omega_{X/O_K}^{Ar} \cdot \omega_{X/O_K}^{Ar}) \geq \frac{\log 2}{6(g-1)}.
\]
\end{Corollary}

\renewcommand{\thesection}{Appendix \Alph{section}}
\renewcommand{\theTheorem}{\Alph{section}.\arabic{Theorem}}
\renewcommand{\theClaim}{\Alph{section}.\arabic{Theorem}.\arabic{Claim}}
\renewcommand{\theequation}{\Alph{section}.\arabic{Theorem}.\arabic{Claim}}
\setcounter{section}{0}

\section{Matrix representation of Laplacian}

In this appendix, we will consider a matrix representation of the Laplacian and
its easy application.

Let $G$ be a metrized graph and $V$ a set of vertices of $G$ such that
$G \setminus V$ is a disjoint union of open segments. Let $E$ be a set of edges
of $G$ by $V$. The length of $e$ in $E$ is denoted by $l(e)$.
Recall that $Q(G, V)$ is a set of continuous functions on $G$ whose
restriction to each edge in $E$ are
quadric polynomial functions, and $M(G, V)$ is
a vector space of measures on $G$ generated by Dirac functions $\delta_v$ at
$v \in V$ and by Lebesgue measures on edges $e \in E$ arising from
the arc-length parameter.
First, we define linear maps $p : Q(G, V) \to \RR^{V}$
and $q : M(G, V) \to \RR^{V}$ in the following ways.
If $f \in Q(G, V)$, then $p(f)$ is the restriction to $V$.
If $\delta_v$ is a Dirac function at $v \in V$, then
\[
q(\delta_v)(v') = \begin{cases}
1 & \text{if $v' = v$} \\
0 & \text{if $v' \not= v$}
\end{cases}
\]
If $dt$ is a Lebesgue measure on a edge $e$ in $E$,
then
\[
q(dt)(v) = \begin{cases}
l(e)/2 & \text{if $v$ is a vertex of $e$} \\
0      & \text{otherwise}
\end{cases}
\]
Next let us define a linear map $L : \RR^{V} \to \RR^{V}$.
For distinct vertices $v, v'$ in $V$,
let $E(v, v')$ be the set of edges in $E$
whose vertices are $v$ and $v'$. Here we set
\[
 a(v, v') = \begin{cases}
   0 & \text{if $E(v, v') = \emptyset $} \\
   {\displaystyle \sum_{e \in E(v, v')} \frac{1}{l(e)}} & \text{otherwise}
 \end{cases}
\]
for $v \not= v'$. Moreover, we set
\[
   a(v, v) = -\sum\begin{Sb} v' \in V \\ v' \not= v \end{Sb} a(v, v').
\]
Let $L : \RR^{V} \to \RR^{V}$ be a linear map
defined by a matrix $(-a(v, v'))_{v, v' \in V}$, i.e.,
if we denote $q(\delta_v)$ by $e_v$, then
$L(e_v) = -\sum_{v' \in V} a(v, v') e_{v'}$.
Thus, we have the following diagram:
\[
\begin{CD}
Q(G, V) @>{\Delta}>> M(G, V) \\
@V{p}VV                        @VV{q}V \\
\RR^{V}      @>>{L}>      \RR^{V}
\end{CD}
\]
Then, we can see the following proposition as remarked in \cite[(5.3)]{BGS}.

\begin{Proposition}
\label{prop:commutativity:L:p:q:Delta}
The above diagram is commutative, i.e., $L \circ p = q \circ \Delta$.
\end{Proposition}

\Proof
Let $f \in Q(G, V)$.
First, let us consider two special cases of $f$.

Case 1 : A case where $f$ is a linear function on each edge in $E$.
By the definition  of $\Delta$, we can see that
\[
\Delta(f) = - \sum_{v \in V}\left(
\sum\begin{Sb} v' \in V \setminus \{ v \} \\
E(v, v') \not= \emptyset \end{Sb}
\left( \sum_{e \in E(v, v')} \frac{f(v') - f(v)}{l(e)} \right)
\right) \delta_v.
\]
On the other hand, by the definition of $a(v, v')$,
\[
\sum\begin{Sb} v' \in V \setminus \{ v \} \\
E(v, v') \not= \emptyset \end{Sb}
\left( \sum_{e \in E(v, v')} \frac{f(v') - f(v)}{l(e)} \right)
= \sum_{v' \in V} a(v, v') f(v').
\]
Therefore, we have
\[
\Delta(f) = \sum_{v \in V} \left(
\sum_{v' \in V} - a(v, v') f(v') \right) \delta_v,
\]
which shows us $q(L(f)) = L(p(f))$.

\medskip
Case 2 : A case where there is $e \in E$ such that $f \equiv 0$ on
$G \setminus e$. Let $v, v'$ be vertices of $e$ and
$\phi : [0, l(e)] \to e$ be the arc-length parameterization of $e$
with $\phi(0) = v$ and $\phi(l(e)) = v'$.
Since $f(v) = f(v') = 0$, $f$ can be written in the form
$f(t) = at(t-l(e))$, where $t$ is the arc-length parameter and
$a$ is a constant. Thus,
\[
\Delta(f) = al(e) \delta_v + al(e) \delta_{v'} - 2a dt.
\]
Therefore, $q(\Delta(f)) = 0$, which means that
$q(\Delta(f)) = L(p(f))$.

\medskip
Let us consider a general case. Let $f_0$
be a continuous function on $G$ such that
$f_0$ is a linear function on each $e \in E$ and
$f_0(v) = f(v)$ for all $v \in V$. Then,
$f - f_0$ can be written by a sum
of functions $f_1, \ldots, f_k$ as in the case 2, i.e.,
\[
f = f_0 + f_1 + \cdots + f_k
\]
and $f_i$ ($1 \leq i \leq k$) is zero on the outside of some edge.
By the previous observation, we know
$q(\Delta(f_i)) = L(p(f_i))$ for all $i = 0, 1, \ldots, k$.
Thus, using linearity of each map, we get our lemma.
\QED

As a corollary, we have the following.

\begin{Corollary}
Let $D = \sum_{v \in V} d_v v$ be a divisor on $G$,
$\mu \in M(G, V)$, and $g \in Q(G, V)$ such that
\[
\int_G \mu = 1 \quad\text{and}\quad
\Delta(g) = \delta_D - (\deg D)\mu.
\]
Then, we have
\[
d_v + \sum_{v' \in V} a(v, v')g(v')  = (\deg D)q(\mu)(v)
\]
for all $v \in V$.
\end{Corollary}

\Proof
Applying $q$ for $\Delta(g) = \delta_D - (\deg D)\mu$ and
using Proposition~\ref{prop:commutativity:L:p:q:Delta},
we have
\[
q(\delta_D) - L(p(g)) = (\deg D)q(\mu).
\]
Thus, by the definition of $L$, we get our corollary.
\QED

\bigskip
Let $k$ be an algebraically closed field,
$X$ a smooth projective surface over $k$, $Y$ a smooth projective curve over
$k$,
and $f : X \to Y$ a generically smooth semi-stable curve of genus $g \geq 1$
over $Y$.
Let $\CV(f)$ be the set of all critical values of $f$ and $y \in \CV(f)$.
Let $G_y$ be the metrized graph of $f^{-1}(y)$ as in
\S\ref{sec:metrized:graph:green:function:admissible:pairing}.
Let $V_y$ be a set of vertices coming from irreducible curves in $f^{-1}(y)$.
For $v \in V_y$, the corresponding irreducible curve is denoted by $C_v$.
Let $K_y$ be the divisor on $G_y$ defined by
$K_y = \sum_{v \in V_y} (\omega_{X/Y} \cdot C_v) v$,
$\mu_y$ the admissible metric of $K_y$, and $g_{\mu_y}$ the Green
function of $\mu_y$.
In this case, the map $L_y : \RR^{V_y} \to \RR^{V_y}$ defined in
the above is given by a matrix $\left(-(C_v \cdot C_{v'})\right)_{v, v' \in
V_y}$.
Thus, the above corollary implies the following proposition.

\begin{Proposition}
\label{prop:admissible:with:fiber}
Let $D$ be an $\RR$-divisor on $X$ and $C_v$ the irreducible curve
in $f^{-1}(y)$ corresponding to $v \in V_y$. Then,
\[
(D \cdot C_v)_a = (D \cdot F) q(\mu_y)(v),
\]
where $F$ is a general fiber of $f$.
In particular, $(D \cdot C_v)_a$ does not depend on the choice
of compactification of $D$.
\end{Proposition}


\end{document}